\documentclass[12pt]{article}
\usepackage{epsfig}
\title{Mostow Rigidity Made Easier}
\author{Richard Evan Schwartz \thanks{Supported by
    N.S.F. Grant DMS-2505281}}

\newtheorem{theorem}{Theorem}[section]

\newtheorem{lemma}[theorem]{Lemma}

\newtheorem{corollary}[theorem]{Corollary}

\def\startproof{{\bf {\medskip}{\noindent}Proof: }}

\def\endproof{$\spadesuit$  \newline}

\def\C{\mbox{\boldmath{$C$}}}%
\def\E{\mbox{\boldmath{$E$}}}%
\def\H{\mbox{\boldmath{$H$}}}%
\def\I{\mbox{\boldmath{$I$}}}%
\def\N{\mbox{\boldmath{$N$}}}%
\def\R{\mbox{\boldmath{$R$}}}%

\begin{document}
\maketitle

\section{The Main Result}

Mostow's rigidity theorem [{\bf M\/}] is one of the most famous
and spectacular results
about hyperbolic manifolds.
The traditional proofs (and also Gromov's proof [{\bf T\/}]) rely on a fair
amount of real analysis. The student who wants to
learn the result {\it all the way to the bottom\/} is in for an arduous journey.
In this article I  give a
proof of Mostow rigidity that is self-contained modulo undergraduate
real analysis.
The proof should be accessible to first-year
graduate students interested in geometry and topology.  My approach
has a lot in common with existing proofs, but it is analytically
lighter.  The ideas here also overlap with [{\bf S\/}].

Given metric spaces $(X_1,d_1)$ and $(X_2,d_2)$, a map
$H: X_1 \to X_2$ is {\it BL\/} (bi-Lipschitz) if $H$ is a bijection
and if there is some
$K \geq 1$ such that
\begin{equation}
K^{-1} d_1(x,y) \leq d_2(H(x),H(y)) \leq K d_1(x,y), \hskip 30 pt \forall x,y \in X_1.
\end{equation}
When $K=1$, the map $H$ is an {\it isometry\/}.
I will prove the following theorem.

\begin{theorem}[Mostow]
  \label{mmm}
If $M_1$ and $M_2$ are 
compact hyperbolic $3$-manifolds and
$f: M_1 \to M_2$ is BL, 
then there is an isometry
$g: M_1 \to M_2$.
So, diffeomorphic compact hyperbolic
$3$-manifolds are isometric.
\end{theorem}
This is a somewhat limited version of Mostow's original Theorem.
In \S \ref{scope}, I discuss
how the same arguments prove the full-blown theorem.
\newline
\newline
{\bf Paper Overview:\/}
The paper is organized in a top-down manner, so that the
big result, Mostow's Theorem, is presented first and proved modulo
three medium-sized results.  The medium-sized results are then
proved modulo some small results, and then finally the
small results are proved.
Here is a chapter-by-chapter account.
\begin{itemize}
\item 
  In \S 1 (this chapter) we introduce background material, and
then prove Mostow Rigidity modulo
Theorem 1.2, Lemma 1.3, and Theorem 1.7.
\item In \S 2 we prove Theorem 1.2 and Lemma 1.3 using classical
  hyperbolic
  geometric arguments, including the Morse Lemma.
  There is nothing new here, but I try to present it well.
  \item In \S 3 we prove Theorem 1.7 modulo several auxiliary results
  that have an analytic flavor.
\item In \S 4 we review some basic measure-theoretic tools that are
  needed for the auxiliary results left over from \S 3.
 \item In \S 5 we use the tools from \S 4 to prove the auxiliary
   results left over from \S 3.
 \end{itemize}

\noindent
{\bf Acknowledgements:\/}
I thank Sujung Jo, Dan Margalit, and ChatGPT for helpful
conversations.
I would also like to thank the anonymous referee for making
helpful expository suggestions.

\subsection{Hyperbolic Geometry Background}
\label{backg}

\noindent
{\bf Hyperbolic Space:\/}
Let $\C$ denote the complex numbers.
{\it Hyperbolic $3$-space\/}, denoted $\H^3$, is a metric
space modeled on the upper
half space in $\C \times (0,\infty)$.
Technically, $\H^3$ is a complete Riemannian manifold
of constant negative sectional curvature.  You don't need
to know this to read the proof of Mostow Rigidity.

Here are some useful facts about the hyperbolic metric. \\[9pt]
{\bf F1.\/} When $a,b \in \C$ and $a \not =0$, the map
$(z,t) \to (az+b,|a|t)$ acts isometrically on $\H^3$\\[9pt]
{\bf F2.\/} The geodesics in $\H^3$, which are the length minimizing paths, are
either vertical rays or semicircles that meet $\C \times \{0\}$
at right angles.\\[9pt]
{\bf F3.\/} Letting $\ell$ and $\ell_{\E}$
denote hyperbolic and Euclidean arc-length respectively, we have
$a\ell \leq \ell_E \leq b \ell$ on $\C \times [a,b]$.\\[9pt]
{\bf F4.\/} The map
$\phi(p)=(0,\|p\|)$ carries $p \in \H^3$ to the point on
$\gamma=\{0\} \times (0,\infty)$ nearest to $p$.
If $t\in (0,1)$, the distance from $(1,t)$ to $\gamma$ is less than $\ln(1/t)+1$.
\newline
\newline
{\bf Hyperbolic Manifolds:\/}
All manifolds are assumed to be compact.
Let $\I$ be the group of isometries of $\H^3$.
A {\it hyperbolic manifold\/} is any quotient of the form
$M=\H^3/\Gamma$, where $\Gamma$ is a subgroup of $\I$.
We only
make this definition when $M$ is actually a manifold.
What this means is that we have a universal covering map
$\pi: \H^3 \to M$, 
and $\Gamma$ is the {\it deck group\/}.  $\Gamma$ is usually
called a {\it co-compact torsion-free lattice\/}, but we will
call $\Gamma$ a {\it nice lattice\/}. The metric on $M$ is such that
$\pi$ is a local isometry.
\newline
\newline
{\bf Conformal Transformations:\/}
Let ${\bf S\/}=\C \cup \infty$ be the Riemann sphere. 
A {\it generalized circle\/} is either a round circle
in $\C$ or else a set of the form $L \cup \infty$
where $L$ is a straight line in $\C$.
A {\it conformal transformation\/} is a homeomorphism
of ${\bf S\/}$ which maps generalized circles to
generalized circles.   When such a map is orientation
preserving, it has the form
\begin{equation}
  \label{mob}
  z \to \frac{az+b}{cz+d}, \hskip 30 pt a,b,c,d \in \C, \hskip 20 pt ad-bc=1.
\end{equation}
As a special case, a {\it homothety\/} is a map of the form
$z \to az+b$ where $a \in (0,\infty)$.    Such maps preserve
the directions of lines in $\C$.

Every isometry of $\H^3$ extends to give a conformal
transformation of ${\bf S\/}$ and every conformal transformation
of ${\bf S\/}$ arises this way.  For this reason, we will sometimes
abuse the notation and speak of a conformal transformation $h$
as being a member of $\I$; technically we are referring here to the
isometry that extends $h$.
\newline
\newline
{\bf Affine Transformations:\/}
An affine map of $\R$ is a map $x \to ax+b$ with $a \not =0$.
Interpreting $\C$ as $\R^2$, we say that a
{\it real affine transformation\/} of $\C$ is a map of the form
$A(v)= T(v)+w$ where $T$ is an invertible real linear transformation
and $w \in \C$.    Most real affine transformations are not conformal.
However, if there are round circles
$C_1,C_2$ such that $A(C_1)=C_2$ then $A$ is in fact conformal.
\newline
\newline
{\bf Equivariance and Lifting:\/}
A BL map $H: \H^3 \to \H^3$ is {\it equivariant\/} if there
is a pair of nice lattices $\Gamma_1,\Gamma_2$ such that
\begin{equation}
  H \Gamma_1 H^{-1}=\Gamma_2.
\end{equation}
Here is another way to express this condition.  For all
$\gamma_1 \in \Gamma_1$, there
exists a $\gamma_2 \in \Gamma_2$ such that
$H \circ \gamma_1=\gamma_2 \circ H$.

Here is how this arises in Mostow Rigidity.
We have $M_j=\H^3/\Gamma_j$ where
$\Gamma_j$ is a nice lattice.
Suppose $f: M_1 \to M_2$ is BL.
Then $f$ has an equivariant
BL {\it lift\/} $H: \H^3 \to \H^3$ which
conjugates $\Gamma_1$ to $\Gamma_2$
and is such that
$$f \circ \pi_1=\pi_2 \circ H.$$
Note that if $H$ is $K$-BL then so is
$fHg$ for any $f,g \in \I$.

We make similar definitions for a
homeomorphism $h: {\bf S\/} \to {\bf S\/}$
(which we usually call a {\it homeo\/})  The
homeo $h$ is
equivariant if $h\Gamma_1h^{-1}=\Gamma_2$.
We will often work with pairs $(H,h)$ where
$h$ is a continuous extension of $H$ to ${\bf S\/}$.
(We prove the existence of such extensions below.)
In this case, we call
the pair $(H,h)$ {\it equivariant\/} if both
maps are simultaneously equivariant with
respect to the same pair of nice lattices.
\newline
\newline
{\bf Tame Sequences:\/}
We call a sequence $\{g_n\}$ in $\I$
{\it tame\/} if 
$\{g_n(p)\}$ is bounded for each
$p \in \H^3$.   In this situation
we can pass to a subsequence and arrange that
there are $4$ points
$p_1,p_2,p_3,p_4 \in \H^3$, say the
vertices of a regular tetrahedron, so that
all $4$ sequences $\{g_n(p_j)\}$ converge.
But then $\{g_n\}$ converges in $\I$. In short,
a tame sequence converges on a subsequence.

Here is a criterion for tameness.
Suppose $a,b,c \in {\bf S\/}$ are $3$ distinct points, and
$\{g_n(a)\}$, 
$\{g_n(b)\}$, 
$\{g_n(c)\}$ all converge to
$3$ distinct points of ${\bf S\/}$.
Then $\{g_n\}$ is tame.
This derives from the
fact that (up to specifying whether the map preserves
or reverses orientation) an element of $\I$ is
determined by what it does on $3$ distinct
points of ${\bf S\/}$. The takeaway here is that
we can tell that a sequence
is tame by looking at how it acts on the ideal boundary.

Below we will generalize these ideas for
sequences of $K$-BL maps.

\newpage

\subsection{A Lesson in Calculus}

Suppose that $h: \R \to \R$ is a function
that is differentiable at some point $a \in \R$.
After translating we normalize so that $a=0$
and $h(a)=0$.  Then what we are saying is
that there is some constant $A$ such that
\begin{equation}
  \lim_{x \to 0} \frac{h(x)}{x}=A.
\end{equation}
It follows that, for any nonzero real number $u$, we have
\begin{equation}
  \lim_{n \to \infty} \frac{h(u/n)}{u/n}=A.
\end{equation}
Put another way
\begin{equation}
  \lim_{n \to \infty} n h(u/n) = Au.
  \end{equation}
  We introduce auxiliary functions
  \begin{equation}
    f_n(x)=nx, \hskip 30 pt
    g_n(x)=x/n.
  \end{equation}
  We define
  \begin{equation}
    \label{SIM}
    h_n=f_n \circ h \circ g_n.
  \end{equation}
  Then we have
  \begin{equation}
    \lim_{n \to \infty} h_n(u)=Au.
  \end{equation}
  This last equation holds for all $u \in \R$, even when $u=0$.
Hence  $h_n \to h'$ where $h'$ is the
  linear function $h'(u)=Au$.
  Given the way that limits work -- and I invite
  you to think about this -- the convergence
  is uniform over compact sets of $\R$.
  In short $h_n \to h'$ uniformly
  on compact sets.  When $A=0$ the
  map $h'$ is the zero map. Otherwise
  $h'$ is a linear map.
  This is how we want to think about the derivative: If
  you \underline{zoom in} to a point of differentiability,
  you get a linear map in the limit.
  
  Before moving on, let us extend the lesson a bit.
  Suppose that $h: \C \to \C$ is a map normalized
  so that $h(0)=0$ and the directional derivative
  of $h$ exists in the direction of a line $L$
  through the origin.  Then we can
  define $h_n$ as above.  If we know that
  the limit $h'$ exists, then we can say that
  $h'|_L$ is a linear
  map from $L$ into $\C$.  The whole map
  $h'$ might be a mystery but the restriction
  to $L$ is not a mystery.
  Our constructions below are all built around this idea.
  
  \subsection{Extensions of Equivariant BL Maps}
  \label{ext}
  
In \S 2 we prove the following well-known result.
  
\begin{theorem}
  \label{extend}
A  BL map $H$ of $\H^3$ extends continuously to a homeo $h$
of ${\bf S\/}$.  If  $H$ is equivariant with respect to a pair
$\Gamma_1,\Gamma_2$ of nice lattices then so is $h$.
\end{theorem}

For the subsequent results in this section we work exclusively
with equivariant pairs $(H,h)$.
We will show in the equivariant case that the map $h$
from Theorem \ref{extend} is a conformal transformation.
This means that $h$ extends to an equivariant isometry
$h^*: \H^3 \to \H^3$ which in turn defines
an isometry $g: M_1 \to M_2$.
So, proving Mostow Rigidity boils down to
proving that $h$ is a conformal transformation.
We will move $h$ around by isometries
until it spills its secret.

We write $(H,h) \sim (H',h')$ if there are
$f,g \in \I$ so that
$H'=fHg$ and $h'=fhg$.
A sequence
$\{(H_n,h_n)\}$ is {\it derived from\/} $(H,h)$
if we have $(H_n,h_n) \sim (H,h)$ for all $n$.
Compare Equation \ref{SIM}.
Generalizing the definition given in \S \ref{backg},
we call a sequence $\{\phi_n\}$ of $K$-BL maps
{\it tame\/} if
$\{\phi_n(p)\}$ is bounded for each $p \in \H^3$. We call
$\{(H_n,h_n)\}$ tame if $\{H_n\}$ is tame.
In \S 2 we prove the following generalization of the
tameness criterion in \S \ref{backg}.

\begin{lemma}
  \label{MORSE}
  Suppose that $\{(H_n,h_n)\}$ is a sequence derived from $(H,h)$ and
  there are $3$ distinct points $a,b,c \in {\bf S\/}$
  such  that $\{h_n(a)\}, \{h_n(b)\}, \{h_n(c)\}$ converge
  to $3$ distinct points in ${\bf S\/}$.
  Then $\{(H_n,h_n)\}$ is tame.
\end{lemma}

Before we implement our calculus lesson we need one more result,
a subtle and powerful one.

\begin{lemma}
  \label{equ}
  Let  $\{(H_n,h_n)\}$ be a tame sequence derived from $(H,h)$.
 Then  $(H,h) \sim (H',h')$ where, on a subsequence,
  $H_n$ converges to  $H'$ uniformly on compact subsets of $\H^3$
  and $h_n$ converges to $h'$ uniformly on ${\bf S\/}$.  
\end{lemma}

\startproof
Let $H_n=f_nHg_n$ and $h_n=f_nhg_n$.
Note that $\{f_n\}$ and $\{g_n\}$ might not be tame.
For any $g \in \Gamma_1$ we have
$f_g \in \Gamma_2$ so that $H=f_gHg$ and
$h=f_ghg$.   Combining these equations with the
fact that $\H^3/\Gamma_1$ is compact,
we can \underline{also} write $H_n=f_n'Hg_n'$ and
$h_n=f_n'hg_n'$,  where $\{g_n'\}$ is tame.
Here
$$f_n'=H_n \circ (g_n')^{-1} \circ H^{-1}.$$
Since $\{H_n\}$ and $\{(g_n')^{-1}\}$ are tame so is
$\{f_n'\}$.   As noted in \S \ref{backg}, we have
$f_n' \to f' \in \I$ and
$g_n' \to g' \in \I$ on a subsequence.
So, we get the desired convergence by setting
$h'=f'hg'$ and $H'=f'Hg'$.
\endproof

\subsection{Zooming In}

Let $(H,h)$ be an equivariant pair.  Our goal is to
show that $h \in \I$.
We normalize so that
$h(\infty)=\infty$ and interpret $h$ as a homeo of $\C$.

\begin{lemma}
  \label{RtoC}
  If $h$ is a real affine transformation
  then $h \in \I$.
\end{lemma}

\startproof
The map $h$ is equivariant with respect to nice
lattices $\Gamma_1$ and $\Gamma_2$.
For any $f \in \Gamma_1$ there is some
$g \in \Gamma_2$ such that
$gh=hf$.  We can choose a line $L$ and an
element $f \in \I$ so that $C_1=f(L)$ is a round circle.
But then $C_2=h(f(L))$ is a (possibly non-circular) ellipse.
On the other hand $h(L)$ is a line and
$C_2=g(h(L))$ is either a line or a round circle.
The only way the two descriptions match is if
$C_2$ is a round circle.
But then $h$ maps the circle $C_1$ to the circle $C_2$.
As we mentioned in
\S \ref{backg}, this forces $h \in \I$.
\endproof

A line $L \subset \C$ is {\it good\/} for $h$ if
the restriction $h|_L$ is an affine map.
This means that $h(L)$ is a line and $h$ scales lengths by a
constant factor on $L$.
We say that the
direction $D$ is {\it good\/} for $h$ if every line of direction $D$
is good for $h$.   
The affine maps here may depend on the line but, since $h$ is a
homeo, the image under $h$ of this parallel line family is another
parallel line family.

\begin{lemma}
  \label{2lines}
  If two directions are good for $h$ then $h \in \I$.
\end{lemma}

\startproof
There are real affine transformations $A,B$ so that
the horizontal and vertical directions for $\phi=AhB$ are
good, and $\phi$ is the identity on the union of the
two coordinate axes.  But $\phi(x,y)$ is on the same
horizontal line as $(0,y)$ and on the same vertical line as
$(x,0)$.  Hence $\phi(x,y)=(x,y)$.  In short, $\phi$ is the identity.
Since $AhB$ is the identity,
$h=A^{-1}B^{-1}$.   Hence $h$ is a real affine map. Hence
$h \in \I$ by Lemma \ref{RtoC}.
\endproof

Choose $z \in \C$.
Specializing the constructions from \S \ref{ext}, we
call $\{(H_n,h_n)\}$ the {\it zoom sequence based on\/} $(h,z)$ if
$H_n=f_nHg_n$ and
$h_n=f_nhg_n$, where
\begin{itemize}
 \item $f_n$ is the homothety that
fixes $h(z)$ and scales distances by $n$. 
 \item $g_n$ is the homothety that
  fixes $z$ and scales distances by $1/n$.
\end{itemize}
We write $h\to h'$ if $\{(H_n,h_n)\}$ is tame, and
$(H_n,h_n)$ converges on a subsequence to $(H',h')$ in the sense
of Lemma \ref{equ}.   For ease of notation we omit mention of $H$,
even though $H$ is lurking in the background at every step.

By Lemma \ref{equ}, we have
$h \sim h'$ if
$h \to h'$. Since $h' \sim h$, the map $h'$ is
also equivariant.
We write
\begin{equation}
  h \Rightarrow h^{(n)} \hskip 20 pt
  {\rm if\/} \hskip 20 pt h \to h' \to h'' \to ... \to h^{(n)}.
  \end{equation}
By several applications of Lemma \ref{equ} we have
$h \sim h^{(n)}$ if $h \Rightarrow h^{(n)}$.  So, if
$h \Rightarrow h^{(n)}$ then $h^{(n)}$ is equivariant.  Also, if
$h \Rightarrow h^{(n)} \in \I$ then $h \in \I$.

The {\it directional derivative\/} of $h$ at $z$ in the direction $v$,
when it exists, is given by
\begin{equation}
  D_vh(z)=\lim_{t \to 0} \frac{h(z+vt)-h(z)}{t}.
\end{equation}
We call this a {\it rational\/} directional derivative if the real and
imaginary
parts of $v$ are rational.
A point $z \in \C$ is
an {\it asterisk\/} for $h$ if every rational directional derivative
exists at $z$ and if $D_1h(z) \not =0$.
In \S 3-5 we prove the following result.

\begin{theorem}
  \label{one}
  $h$ has an asterisk.
      \end{theorem}

      \begin{corollary}
  \label{star1}
  Let $D$ be a direction.  Then $h \Rightarrow h''$
  where $D$ is good for $h''$.
\end{corollary}

\startproof
Let $z$ be an asterisk for $h$. Let $\{h_n\}$ be the zoom sequence
based on $(h,z)$.  Let $\cal Q$ denote the set of lines of
rational slope through $z$ and let $\cal R$ denote the set of all
lines through $z$.
 Remembering the calculus
lesson, we see that
the restriction of $h_n$ to each line of $\cal Q$
converges to an affine map.
In particular $h_n(z+m) \to h(z)+mD_1h(z)$ for $m=0,1,2$,
and these points are distinct.
So $\{h_n\}$ is tame by Lemma \ref{MORSE}.  Hence
$h \to h'$ and the restriction of the homeo $h'$ to every line of $\cal Q$ is
a nontrivial affine map.  Hence every line of $\cal Q$
is good for $h'$.
Since $h'$ is a homeo, and $\cal Q$ is dense in $\cal R$,
 all lines of $\cal R$ are good for $h'$.

Choose $z' \in \C$ so that the line through $z$ and $z'$, which we
denote as $\overline{zz'}$, has direction $D$.
Let $h'_n=f_n'h'g_n'$.   The sequence
$\{h'_n\}$ is the zoom sequence based on $(h',z')$.
The restriction of
$h'_n$ to $\overline{zz'}$ is independent of $n$.
So, by Lemma \ref{MORSE},
$\{h'_n\}$ is tame.  Every line in
$(g'_n)^{-1}({\cal R\/})$ is good for $h_n'$,
and as $n \to \infty$ this set converges to the
set of lines in direction $D$. (The point $(g_n')^{-1}(z)$ common to
all these lines moves to $\infty$ along the line $\overline{zz'}$.)
Hence $D$ is good for $h''$.
\endproof

If $D$ is good for $h$ and $h \Rightarrow h''$
then $D$ is also good for $h''$.
So, applying Corollary \ref{star1} twice, for different
directions, we get $h \Rightarrow h^{(4)}$ where
two directions are good for $h^{(4)}$, and $h^{(4)}$ is equivariant.
But then $h^{(4)} \in \I$ by Lemma \ref{2lines}.
By Lemma
\ref{equ} again, $h \in \I$.
Our proof of Mostow Rigidity is done.

\subsection{Discussion}

The idea of the asterisk is the most unusual ingredient
in our proof.  I cannot remember how I thought of it, but
let me explain how it functions in the proof, and what
considerations might lead to it.  The
essential idea behind the above kind of proof of Mostow
Rigidity is that zooming in to a point that has some regularity
reveals some extra structure of the map.  Let's discuss
what the phrase {\it some regularity\/} might mean.

At one extreme, {\it some regularity\/} might mean
that the partial derivatives exist at the point.  If the
partial derivatives are zero, this information is
completely useless to us.  If the partial derivatives
are nonzero, then when we zoom in and take a suitable
limit, we find (after conjugating) that our map is
actually a linear transformation when restricted to
two lines.  This is nice to know, but does not move us
much towards a proof of Mostow Rigidity.
On the positive side, this
amount of regularity is not too hard to establish.

At the other extreme, {\it some regularity\/}
might mean that the map is actually differentiable
at a point.  As is well known, the map $h$ from
Theorem \ref{extend} is {\it quasi-conformal\/}.
Such maps are known to be differentiable almost everywhere,
and the derivative is non-singular almost everywhere.
See  [{\bf LV\/}].  Given a
point of nonsingular differentiability at $z$, we get a swifter
proof. We consider the zoom sequence $\{h_n\}$ based
on $(h,z)$ and we get $h \to h'$ where $h'$ is real affine
and equivariant -- and hence in $\I$.
The difficulty with this alternate proof is that the
analysis behind the differentiability result just
quoted is (to me) rather formidable.

The asterisk idea is a happy compromise between
these two extremes. On the one hand, the amount
of work needed to show the existence of an asterisk
is not much more than the
amount needed to show that the partial
derivatives exist and are nonzero at a point.
On the other hand, when we zoom into an
asterisk we get control over what happens in
all rational directions.  This is an intermediate amount
of information, compared to the two extremes,
but it is enough to push through the proof.
The rest of this section will give some
variations on a theme, showing other
ways to convert the asterisk idea into a proof.
\newline
\newline
\noindent
{\bf Three Fair Directions:\/}
Let us say that $h$ is {\it fair\/} on a line $L$ if
$h(L)$ is a straight line.   We say that a direction
$D$ is {\it fair\/} for $h$ if $h$ is fair for every
line in the direction $D$.
Corollary \ref{star1} says
in particular that for any direction $D$ we have
$h \Rightarrow h''$, where $D$ is fair for $h''$.
Applying Corollary \ref{star1} three times we
get $h \Rightarrow h^{(6)}$, where $3$ different
directions are fair for $h^{(6)}$.  It is a fun exercise
to show that this forces the equivariant $h^{(6)}$ to be real affine
-- and hence in $\I$.
\newline
\newline
{\bf Zoom Ad Infinitum:\/}
We could apply Corollary \ref{star1} an infinite
number of times and then take a limit
(guaranteed by a result much like Lemma \ref{equ})
and conclude that $h \sim h^*$ where $h^*$ is
fair on every line.  But then the equivariance
implies that $h$ preserves generalized circles and hence is in $\I$.
\newline
\newline
{\bf The Minimal Proof:\/}
The first half of the proof of Corollary \ref{star1} shows that
$h \to h'$ where all lines through a point $z \in \C$
are good for $h'$.  Since $h \sim h'$ we see that
there are two points $a,b \in {\bf S\/}$ such that
the restriction of $h$ to each generalized circle through $a,b$ is a
conformal map.   Since $h$ is equivariant, there is a second
pair $a',b'$ of points, both distinct from $a$ and $b$, such
that the restriction of $h$ to each generalized circle through $a',b'$ is
a conformal transformation.  It is a fun exercise to show that
these conditions force $h \in \I$.  This proof is a bit harder to
generalize to higher dimensions.
\newline
\newline
{\bf Challenge:\/} I was hoping to figure out a proof which combines
equivariance with a differentiability result even weaker than
Theorem \ref{one} but couldn't get this to work. Can you?

\newpage

\section{The Extension}

\subsection{Proofs modulo the Morse Lemma}

In this chapter we prove
Theorem \ref{extend} and Lemma \ref{MORSE}.
In this first section we reduce these
results to a well-known geometric result
called the Morse Lemma.  See [{\bf T\/}, Prop. 5.9.2].
We work with the $K$-BL map
$H: \H^3 \to \H^3.$

Given a geodesic $\gamma \subset \H^3$ let
$N_r(\gamma)$ denote the set of points
in $\H^3$ at most $r$ from $\gamma$.
The set $N_r(\gamma)$ is called the
$r$-{\it tubular neighborhood\/} of $\gamma$.

\begin{lemma}[Morse]
  Let $\gamma$ be a geodesic.  Then
  there is a unique geodesic
  $\gamma'$ such that $H(\gamma) \subset N_{K'}(\gamma')$.
  The constant $K'$ only depends on $K$.
\end{lemma}

 Our proof gets
the constant $K'=4K^3+2K+1$, but this is not optimal.
The choice of $K'$, however,
has no impact on the proof of Theorem \ref{extend}.
The Morse Lemma holds in great generality,
with good control on $K'$. See [{\bf GS\/}].
\newline
\newline
{\bf Proof of Theorem \ref{extend}:\/}
The Morse Lemma induces a map
of the set of oriented geodesics.  If $\gamma$ is an oriented
geodesic
then $\gamma \to \gamma'$ where
$\gamma'$ is the unique geodesic such that
$H(\gamma) \subset N_{K'}(\gamma')$.
The orientation of $\gamma$ combines with $H$ to
determine the orientation of $\gamma'$.
If $\tau$ is another geodesic and
$\tau \to \tau'$ and
$\gamma,\tau$ have a common endpoint,
then $H(\gamma)$ and $H(\tau)$ are asymptotic
to each other in the corresponding direction.
Hence $\gamma'$ and $\tau'$ also have
a common endpoint.  This property lets us define
$h$ as the unique map of the Riemann sphere
${\bf S\/}$ with the
following property: if $\gamma \to \gamma'$ and
$\gamma$ connects
$p$ and $q$, then $\gamma'$
connects $h(p)$ to $h(q)$.
If we apply the construction to $H^{-1}$ we get $h^{-1}$.

It remains to show that $h$ is a homeo.
By construction, $h$ is a bijection.
If $h$ is not a homeo, 
there are $H$-corresponding sequences of
geodesics $\{\gamma_{n}\} \leftrightarrow
\{\gamma_{n}'\}$
such that the endpoints of the geodesics in the one sequence
come together on ${\bf S\/}$ and the endpoints of the geodesics in the other
sequence do not.  But then the one sequence exits every
compact subset of $\H^3$ and the other sequence does not.
This contradicts the BL nature of $H$.
Given the way that the extension only depends on $H$ we see
that $h$ is equivariant if $H$ is equivariant.
\endproof

Before proving Lemma \ref{MORSE} we
need one preliminary lemma.  An
{\it ideal triangle\/} is a union of
$3$ geodesics defined by $3$ distinct
points on ${\bf S\/}$.   Any two
ideal triangles are equivalent under the
action of $\I$.

\begin{lemma}
  \label{bounded}
  Let $\Delta$ be an ideal geodesic triangle in $\H^3$.
  Then for any $R$, the set of points in $\H^3$ within
  $R$ of all three geodesics of $\Delta$ is compact.
\end{lemma}

\startproof
Let $\gamma_{a,b}$ be the geodesic in $\H^3$
whose endpoints are $a,b \in {\bf S\/}$.
We can normalize so that $\Delta$ has vertices
$0,1,\infty$.   Let $p=(z,t) \in \H^3$.
If $t \to 0$ then the distance from $p$ to
one of $\gamma_{0,\infty}$ or $\gamma_{1,\infty}$
tends to $\infty$.  If $t \to \infty$ or
$|z| \to \infty$ then the distance from
$p$ to $\gamma_{0,1}$ tends to $\infty$.
So, $p$ is close to all these geodesics when
$t \sim 1$ and $|z|$ is fairly small.
\endproof

\noindent
{\bf Proof of Lemma \ref{MORSE}:\/}
We can normalize by a pair of elements of
$\I$ so that
$\{(H_n,h_n)\}$ is such that
$h_n(0) \to 0$ and $h_n(1) \to 1$ and $h_n(\infty) \to \infty$.
Further composing with convergent sequences in
$\I$ we can assume that $h_n$ fixes each of $0,1,\infty$.
Let $p \in \H^3$ be arbitrary.   By the Morse Lemma,
$H_n(p)$ is uniformly close to all
three of the geodesics of the ideal triangle with vertices
$0,1,\infty$.  But then, by Lemma \ref{bounded},
$\{H_n(p)\}$ remains within a bounded subset
of $\H^3$.  Hence $\{(H_n,h_n)\}$ is tame.
\endproof

\subsection{The Morse Lemma Modulo the Tube Lemma}

We keep the notation from the previous section.
Let $d$, $\ell$ respectively denote hyperbolic
distance, hyperbolic arc length.
As in {\bf F4\/} from \S \ref{backg}, we
define
\begin{equation}
  \label{PROJMAP}
  \gamma=\{0\} \times (0,\infty), \hskip 30 pt
  \phi(p)=(0,\|p\|) \in \gamma.
\end{equation}
Again, $\phi$ is the map such that
$\phi(p)$ is the point of $\gamma$ closest to $p$.

In the next section we prove the following result.

\begin{lemma}[Tube]
  \label{tube}
  If $r>1$ then $\ell \circ \phi  \leq  e^{-r+1} \ell$ on
  $\H^3-N_r(\gamma)$.
\end{lemma}
Now we use the Tube Lemma to prove the Morse Lemma.
The proof comes in two steps.  The first step is showing
that an arc of $H(\gamma)$ cannot wander too far from
the geodesic connecting its endpoints.   If the arc does
wander far away,
we can use $\phi$ to produce a much shorter
competing path which connects the endpoints.
This will give a contradiction.
After we have this result about finite arcs of $H(\gamma)$, the
second step involves considering a sequence of increasingly
long arcs of $H(\gamma)$ and taking a suitable limit.
The key point is that the geodesics connecting the
endpoints of these arcs all have to line up.

\begin{lemma}
  \label{quasi}
  Let $C=4K^3+2K$.
  Let $\alpha$ be a geodesic segment in $\H^3$.  Then
  $H(\alpha) \subset N_{C}(\gamma')$ for
the geodesic $\gamma'$ through the
endpoints of $H(\alpha)$.
\end{lemma}

\startproof
We normalize so that $\gamma'=\{0\} \times (0,\infty)$, as in the
Tube Lemma.  We show a schematic picture of the construction.

\begin{center}
\resizebox{!}{1.5in}{\includegraphics{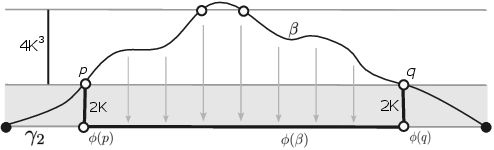}}
\end{center}

Let $\beta=H(\alpha)$.
Suppose $\beta \not \subset N_{C}(\gamma')$.
Then there are $p,q \in \partial N_{2K}(\gamma')$ and an
arc $\beta_{pq}$ of $\beta$ joining $p$ to $q$ that remains outside
$N_{2K}(\gamma')$ and exits $N_{C}(\gamma')$.  Note that
$\ell(\beta_{pq}) \geq 8K^3$.
The path $\phi(\beta_{pq})$ connects $\phi(p)$ to $\phi(q)$.
By the Tube Lemma
\begin{equation}
  \label{TUBE}
  d(\phi(p),\phi(q)) \leq \ell(\phi(\beta_{pq})) \leq \ell(\beta_{pq})
  e^{-2K+1}.
  \end{equation}
By {\bf F4\/} from \S \ref{backg}, we have
$d(p,\phi(p))=d(q,\phi(q))=2K$.  By
Equation \ref{TUBE} and the triangle inequality
\begin{equation}
  \label{morse00}
  d(p,q) \leq d(p,\phi(p))+d(\phi(p),\phi(q))+d(\phi(q),q)  \leq 
  e^{-2K+1}\ell(\beta_{pq})+4K.
\end{equation}

Since $\beta_{pq}$ is the image of a geodesic segment
under a $K$-BL map,
\begin{equation}
  \label{morse01}
  \ell(\beta_{pq})  \leq K \ell(H^{-1}(\beta_{pq})) = K d(H^{-1}(p),H^{-1}(q))
  \leq K^2\, d(p,q).
\end{equation}
Combining Equations \ref{morse00} and \ref{morse01}, and using
$\ell(\beta_{pq}) \geq 8K^3$, we have
\begin{equation}
  \label{morse1}
  \bigg(e^{-2K+1} + \frac{1}{2K^2}\bigg) \ell(\beta_{pq}) \geq
  e^{-2K+1} \ell(\beta_{pq}) + 4K \geq d(p,q) \geq  \frac{\ell(\beta_{pq})}{K^2}.
\end{equation}
Multiplying through by $2K^2/\ell(\beta_{pq})$ and subtracting $1$, we get
$$2K^2 e^{-2K+1} \geq 1.$$
This is false for all $K \geq 1$, a contradiction.
\endproof

\noindent
{\bf Proof of the Morse Lemma:\/}
Let $C=4K^3+2K$ be as in Lemma \ref{quasi}.
We fix an origin $O$ of $\gamma$ and let $\alpha_{n} \subset
\gamma$ be a
geodesic segment of length at least $2n$ centered at $O$.
Let the endpoints be $a_n$ and $b_n$.
Each ray of $\gamma$ emanating from $O$ exits
every compact subset of $\H^3$, so we can choose our
segments so that
$\{H(a_n)\}$ and $\{H(b_n)\}$ both converge to
points $A,B \in {\bf S\/}$.   Let
$\gamma'_n$ be the geodesic through the
endpoints of $H(a_n)$ and $H(b_n)$.
Given our description of
hyperbolic geodesics, and
the fact (from Lemma \ref{quasi}) that $\gamma'_{n}$ comes
 within $C$ of $H(O)$ for all $n$, we cannot have $A=B$.

We normalize so that $A=(-1,0)$ and $B=(1,0)$.
Let $\Delta_r \subset \H^3$ be the ball of radius $r$ about $O$.
Let $\gamma'$ be the geodesic connecting $(-1,0)$ and $(1,0)$.
The endpoints of $\gamma'_{n}$ converge
to the endpoints of $\gamma'$.  By compactness and
Lemma \ref{quasi},
we have the following for sufficiently large $n$:
$$H(\gamma) \cap \Delta_r = H(\alpha_{n}) \cap \Delta_r \subset
N_C(\gamma'_{n}) \cap \Delta_r \subset N_{C+1}(\gamma') \cap \Delta_r
\subset N_{C+1}(\gamma').$$
Letting $r \to \infty$ we see that $H(\gamma) \subset
N_{C+1}(\gamma')$. Now we set $K'=C+1$.
\endproof

\subsection{Proof of the Tube Lemma}

Let $\phi$ and $\gamma$ be as in Equation \ref{PROJMAP}.
Let $\ell_{\E}$ denote Euclidean arc length.   Again,
$\ell$ denotes hyperbolic arc length.
We use the facts {\bf F1\/} -- {\bf F4\/} from \S \ref{backg}.
Looking at the formula for $\phi$ we can see that
$  \ell_{\E}\circ \phi \leq \ell_{\E}$.
That is, $\phi$ does not increase Euclidean
arc length.
Define
\begin{equation}
  \tau=\{1\} \times (0,\infty), \hskip 30 pt
U=\big(\C \times (0,e^{-r+1})\big) \cap \phi^{-1}(\C \times
(1,\infty))
\end{equation}
The set $U$ is open and $\phi(U) \subset \C \times (1,\infty)$.
By {\bf F3\/}, we have $\ell \leq \ell_{\E}$ on $\phi(U)$ and
$\ell_{\E} \leq e^{-r+1}\ell$ on $U$. Hence, on $U$, we have
$$ \ell \circ \phi \leq \ell_{\E} \circ \phi \leq \ell_{\E} \leq e^{-r+1}
\ell.$$

Combining {\bf F3\/} and {\bf F4\/} we see that
$(1,t)$ is less than $r$ units from $\gamma$ when
$t \in [e^{-r+1},\infty)$. Therefore
\begin{equation}
  \label{chop}
  \tau-N_r(\gamma) \subset U.
\end{equation}

Let $\Gamma \subset \I$ denote the stabilizer
subgroup of $\gamma$.  Elements of $\Gamma$ have
the form $(z,t) \to (az,|a|t)$ with $a \in \C-\{0\}$.
Every point of $\H^3 - N_r(\gamma)$
can be mapped into $\tau-N_r(\gamma)$ by an element of
$\Gamma$.  Hence $\H^3-N_r(\gamma) \subset \Gamma(U)$, the
orbit of $U$ under $\Gamma$.
But elements of $\Gamma$ commute with $\phi$. For this reason,
our inequality $\ell \circ \phi \leq e^{-r+1} \ell$ on $U$ holds on all
of $\H^3 - N_r(\gamma)$.
This proves the Tube Lemma.

\newpage

\section{Existence of Asterisks}

\subsection{Reduction to a Technical Lemma}

A set $S \subset \C$ is {\it null\/} if, for every
$\epsilon>0$, there is a countable collection $\{D_j\}$ of disks
such that $S \subset \bigcup D_j$ and
$\sum {\rm area\/}(D_j)<\epsilon$.  Otherwise, $S$ is {\it fat\/}.
For someone who knows about measure theory, we remark that
a null set is one with zero outer measure and a fat set
is one with positive outer measure.

Let $(H,h)$ be an equivariant pair, normalized
so that $h(\infty)=\infty$. The purpose of this chapter
is to prove that $h$ has an asterisk, a point where
all the rational directional derivatives exist and
are nonzero.

\begin{lemma}
  \label{diff}
  The subset of $[0,1]^2$ where $\partial h/\partial x$ does not exist
  is null.   The subset of $[0,1]^2$ where $\partial h/\partial x \not = 0$ is
  fat.
\end{lemma}

\noindent
{\bf Proof of Theorem \ref{one}:\/}
  Tile $\C$ by unit squares.   By symmetry and Lemma \ref{diff}
  the derivative $\partial h/\partial x$ exists in each square outside
  of a null set.  The countable union of null sets is null, so
  $\partial h/\partial x$ exists outside of a null set $S_0 \subset \C$.  By
  rotational symmetry $D_rh$ exists outside a null set for all
  rational vectors.  The countable union
    of these null sets is null, and outside $\bigcup S_r$ all the
  rational directional derivatives of $h$ exist.
A subset of a null set is null, so the
fat set where $\partial h/\partial x \not =0$
intersects the set where all the rational directional
derivatives of $h$ exist.  This gives us an asterisk,
and in fact many.
\endproof

\subsection{Analytic Preliminaries}
\label{BBB}

The rest of the chapter is devoted to the proof of
Lemma \ref{diff}.    We first gather together some analytic
results.

The {\it Borel $\sigma$-algebra\/} is the smallest
collection of subsets of $[0,1]^d$ that contains
all closed subsets and is closed under the
operations of taking
complements, countable intersections,
and countable unions. We only care
about the cases $d=1,2$.
A {\it Borel set\/} is a member of
the Borel $\sigma$-algebra.

\begin{lemma}
  \label{nice}
  Let $\phi: [0,1]^2 \to \R$ be a continuous function.
  Then the set of points in $(0,1)^2$
  where $\partial \phi/\partial x$ exists is
  a Borel set.
\end{lemma}
We give a self-contained proof of this result in \S \ref{niceproof}.
It is a special case of an extremely general
result, [{\bf Z\/}, Prop. 3.3].

We identify the set of horizontal lines in $[0,1]^2$ with the
interval $[0,1]$ in the obvious way. 
\begin{theorem}
  \label{null2}
  Let $S \subset [0,1]^2$ be a set.
    Let  $F_S \subset [0,1]$ be the set
  of horizontal lines $L$ such that
  $S \cap L$ is fat.  If $S$ is null then $F_S$ is null.
  If $S$ is a Borel set and $F_S$ is null then
  $S$ is null.
\end{theorem}
This result is a special case of the Fubini-Tonelli
Theorem [{\bf F\/}, p. 65].
We give a self-contained proof in \S \ref{fullproof}.

Let $\cal J$ denote the set of closed intervals of $[0,1]$.
Suppose $A: {\cal J\/} \to (0,\infty)$ has the property that
$A([0,1]) \geq \sum_j A(J_j)$ when $\{J_j\}$ is a set of
disjoint intervals in $[0,1]$.
An interval $J \subset [0,1]$ is
$N$-{\it stretched\/} if  $|A(J)| \geq N|J|$.
A point $p \in [0,1]$ is {\it stretchy\/} if
for any $N$ there is an interval $J$ centered at $p$
that is $N$-stretched.   Otherwise we call $p$ {\it stiff\/}.

\begin{theorem}
  \label{stiff}
  The set of stretchy points in $[0,1]$ relative to $A$ is null.
\end{theorem}
This result is a bit too idiosyncratic to have a direct proof in
the literature, but it is a quick consequence of  Lemma
\ref{cover1} below, which is a special case of the kind of covering
result used in Folland's proof of the Maximal Theorem. See
 [{\bf F\/}, \S 3.4].
 We give a self-contained proof of Theorem \ref{stiff}
 in \S \ref{stiffproof}.

  Now we come to the crucial analytic definition, the
  definition of an absolutely continuous function.
  \newline
  \newline
  {\bf Definition (AC Functions):\/}
  Suppose $f: [0,1] \to \R$ is continuous.
Let $I=\{I_1,...,I_n\}$ denote a finite list of
intervals of $[0,1]$ having pairwise disjoint
interiors.  We call $I$ a
{\it partial partition\/}.  Let $|I|=\sum |I_k|$.
We define $I_k'$ to be the interval bounded
by the two points of $f(\partial I_k)$.
We define $I'=\{I_1',...,I_n'\}$ and
$|I'|=\sum |I_k'|$.
The function $f$ is {\it AC\/} (absolutely continuous)
if, for
each $\epsilon>0$, there is some
$\delta>0$ such that
$|I|<\delta$ implies that $|I'|<\epsilon$.
\newline

\begin{theorem}
\label{AC}
Suppose $f: [0,1] \to \R$ is AC.  Then
\begin{enumerate}
  \item 
$f$ is differentiable on the complement
of a null set.
\item 
If $f(0) \not = f(1)$ then
$f' \not =0$ on a fat set.
\end{enumerate}
\end{theorem}
Theorem \ref{AC} is an immediate consequence of
what Folland calls {\it The Fundamental Theorem of Calculus
  for Lebesgue Integrals\/}.  See [{\bf F\/}, p. 102].
We give a self-contained proof in \S \ref{ACproof}.

\subsection{Proof of Lemma \ref{diff}}
\label{outerdif}

We work with closed intervals and disks, though sometimes consider
their interiors.
Here is a well-known property of our homeo $h$.
See [{\bf T\/}, \S 5.9].

\begin{theorem}[Disk]
There is a constant $K$,
depending only on $h$,
with the following property.
Let $\Delta\subset \R^2$ be a disk. Then
there are disks $D_1,D_2$ 
so that $D_1 \subset h(\Delta) \subset D_2$ and
${\rm diam\/}(D_2)/{\rm diam\/}(D_1)<K$.
\end{theorem}

\startproof
Suppose $\{\Delta_n\}$ is a sequence where
the best ratio for $h(\Delta_n)$ tends to $\infty$.
Let $\Delta$ be the unit disk.
Composing with homotheties, and passing to a subsequence
we get a derived
sequence $\{(H_n,h_n)\}$ such that
$h_n(0)=0$ and $\{h_n(1)\}$ converges.  The best
ratio for $h_n(\Delta)$ tends to $\infty$.
This sequence is tame by
Lemma \ref{MORSE}.  So, passing to a subsequence again,
we arrange that $h_n$ converges uniformly to a homeo $h'$.
Since $h'$ is a homeo, there are concentric disks
$D_1, D_2$ with $D_1^o \subset h'(\Delta) \subset D_2^o$.
Here $D_j^o$ is the interior of $D_j$.  But then
$D_1 \subset h_n(\Delta) \subset D_2$ for large $n$.
This is a contradiction.
\endproof

Now we come to the key geometric idea.
I learned this idea from [{\bf LV\/}].
When $S \subset \C$ we define
$\alpha(S)$ to be the supremum of all finite sums
$\sum {\rm area\/}(D_i)$ where
$\{D_i\}$ is a collection of disjoint disks contained in $S$.
For each interval $J \subset [0,1]$ let
\begin{equation}
  \label{lv}
  A(J)=\alpha(h([0,1] \times J)).
  \end{equation}
By construction, $A$ satisfies the
hypothesis of Theorem \ref{stiff}.

\begin{theorem}
  \label{heart}
  Suppose $y \in (0,1)$ is a stiff point for $A$.  Let
$\pi: \C \to \R$ be any linear projection.
Then $\pi \circ h$ is AC on $L_y=[0,1] \times \{y\}$.
\end{theorem}

\noindent
{\bf Proof of Lemma \ref{diff}:\/}
Let $h_1={\rm Re\/}(h)$ and
$h_2={\rm Im\/}(h)$.
Let $D_j \subset (0,1)^2$ be the set where $\partial h_j/\partial x$ exists.
By Lemma \ref{nice},
the set $D_j$ is a Borel set.
Let $L_y$ be some horizontal line where $y$ is a stiff point.
By Theorems \ref{heart}
and \ref{AC}, we see that $D_j^c \cap L_y$ is null.
So, $D_j^c \cap L_y$ is null unless $y$ is stretchy, and
the stretchy set is null by Theorem \ref{stiff}.  By Theorem
\ref{null2}, $D_j^c$ is null.
Hence $\partial h/\partial x$ exists outside of the null
set $D_1^c \cup D_2^c$.
Let $T \subset [0,1]^2$ be the set where
$\partial h/\partial x$ exists and is nonzero.
There is a fat set $F \subset [0,1]$ and a linear
projection $\pi$ with the following property.
For each $y \in F$, the map $\pi \circ h$ does not
identify the endpoints of
$[0,1] \times y$.   By Theorem \ref{AC},
the set $T \cap L_y$ is a fat subset of $L_y$ for any $y \in F$.
By Theorem \ref{null2}, the set $T$ is fat.
\endproof
\newline
\newline
{\bf Proof of Theorem \ref{heart}:\/}
Let $Q=[0,1]^2$. 
Let $f=\pi \circ h|_{L_y}$.  
If $f$ is not AC,
we can scale so that there
is a sequence of partial partitions $\{I^n\}$ with
$|I^n|<1/n$ and $|(I^n)'| \geq K$, the constant
from the Disk Theorem.  We can subdivide so that
the intervals in each partition have the same size up
to a factor of $2$.

 For any set $Y$ let $Y^*=h(Y)$.
Fix $n$ and let
$I^n=\{I_1,...,I_k\}$.
Let $\epsilon=\max |I_j|$.  Since $\epsilon<2\min |I_j|$, we have
$k\epsilon<2|I|$. Hence $\epsilon<2/(kn)$.
Let $Q_{\epsilon}=[0,1] \times [y-\epsilon,y+\epsilon]$.
Since $y$ is stiff,
there is some $\Omega$ such that
$  \alpha(Q_{\epsilon}^*)<\Omega \epsilon$.
Let $\Delta_j$ be the open disk having $I_j$ as a diameter.
These disks are disjoint.  The figure shows the situation when $k=3$.

\begin{center}
\resizebox{!}{.95in}{\includegraphics{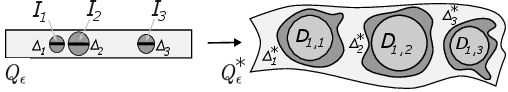}}
\end{center}

Since $h$ is a homeo,
the image sets $\{\Delta_j^*\} \subset Q_{\epsilon}^*$ are also
disjoint.  Hence
\begin{equation}
\label{small0}
\sum_{j=1}^k  \alpha(\Delta_j^*) 
\leq \Omega\epsilon.
\end{equation}
Let $D_{1,j} \subset \Delta_j^* \subset D_{2,j}$ be as in the Disk
Theorem.
We have
\begin{equation}
\label{big2}
\sum_{j=1}^k {\rm diam\/}(D_{2,j}) \geq K,
\hskip 15 pt
\sum_{j=1}^k {\rm diam\/}(D_{1,j}) \geq 1, \hskip 15 pt
\sum_{j=1}^k {\rm area\/}(D_{1,j}) \geq \frac{1}{4k}.
\end{equation}
The first equation, which comes from $|I'| \geq K$, 
implies the second equation, and the second equation
implies the third.  Hence

\begin{equation}
\label{big3}
\frac{2\Omega}{kn}= \Omega \epsilon>\sum_{j=1}^k \alpha(\Delta_j^*) \geq
\sum_{j=1}^k {\rm area\/}(D_{1,j}) \geq \frac{1}{4k}.
\end{equation}
So $n<8\Omega$, a contradiction when $n$ is large enough.
\endproof

\newpage

\section{A Primer on Measure Theory}

\subsection{Measurability}

This chapter proves some standard measure-theoretic results.
The material in this section can be found in [{\bf F\/}, \S 1].

We work in $[0,1]^d$. We only care about $d=1,2$.
A {\it dyadic cube\/} in $[0,1]^d$ is a cube we
get by starting with $[0,1]^d$ and recursively
subdividing cubes into $2^d$ equal-sized
sub-cubes and selecting one of them.
Two dyadic cubes are either nested or else have disjoint interiors.
Say that a {\it carpet\/} is a countable union of dyadic cubes
having pairwise disjoint interiors.
When $S \subset [0,1]^d$ we define the {\it outer measure\/}
\begin{equation}
  \mu(S)=\inf_{S \subset \bigcup Q_n}
  \sum {\rm volume\/}(Q_j)
\end{equation}
The infimum is taken over all carpets containing $S$.

A subset $S \subset [0,1]^d$ is {\it open\/} iff the complement
$S^c=[0,1]^d-S$ is closed.
Every open subset of $[0,1]^d$  is the countable union of dyadic
cubes.  In particular, this is true for open balls.
Using this fact it is an easy exercise to show that
$S$ is null iff $\mu(S)=0$,
and fat iff $\mu(S)>0$.
We also have:
\begin{itemize}
\item  {\it monotonicity:\/} $\mu(A) \leq \mu(B)$ for any  $A \subset
  B$.
    \item {\it subadditivity\/} If $S=\bigcup T_n$ then $\mu(S) \leq
    \sum \mu(T_n)$.
\end{itemize}

A subset $S \subset [0,1]^d$ is {\it measurable\/} if
\begin{equation}
  \label{meas}
  \mu(E \cap S)+\mu(E \cap S^c)=\mu(E)
\end{equation}
for all subsets $E \subset [0,1]^d$.
Note that $S$ is measurable iff $S^c$ is measurable.
It follows directly from monotonicity and subadditivity
that null sets are measurable.

We defined Borel sets in \S \ref{BBB}.

\begin{theorem}
  \label{borel}
  \label{sandwich}
Let $S$ be a Borel set.  Then the following is true.
  \begin{enumerate}
  \item $S$ is measurable.
  \item If $S=\bigcup S_n$, an increasing union of Borel sets, then
         $\mu(S)=\lim \mu(S_n)$. 
  \item For any $\epsilon>0$ we have a compact set $K$ and an open
    set $U$ such that $K \subset S \subset U$ and $\mu(U-K)<\epsilon$.
    \end{enumerate}
  \end{theorem}

  We prove this result through a series of lemmas.

  \begin{lemma}
    \label{cube}
    Dyadic cubes are measurable.
  \end{lemma}

  \startproof
  Let $A$ be a dyadic cube.  Let $E \subset [0,1]^d$ be an arbitrary
  set.  Any carpet containing $E$ can be further subdivided so that
  each of its cubes is contained in either $A$ or $A^c$.
  The union of the former gives a sub-carpet covering
  $E \cap A$ and the union of the latter gives a sub-carpet covering
  $E \cap A^c$.  Since this holds for all covers of
  $E$ we have
$\mu(E)=\mu(E \cap A) + \mu(E \cap A^c)$.
\endproof

\begin{lemma}
  If $A, B$ are measurable then so are $A \cap B$ and $A \cup B$
  and $A-B$.
  \end{lemma}

  \startproof 
  Note that $A^c$ and $B^c$ are also measurable.
  Since $(A \cup B)^c=A^c \cap B^c$ and
  $A-B=A \cap B^c$, it suffices to prove that $A \cap B$ is
  measurable.   Let $R_k$ be the region shown in the Venn
    diagram for $A,B,E$ as in the figure.  Let
    $[k_1...k_m]=\mu(R_{k_1})+...+\mu(R_{k_m})$.
We want to prove
    $[4567]=[456]+[7]$.

\begin{center}
\resizebox{!}{1.8in}{\includegraphics{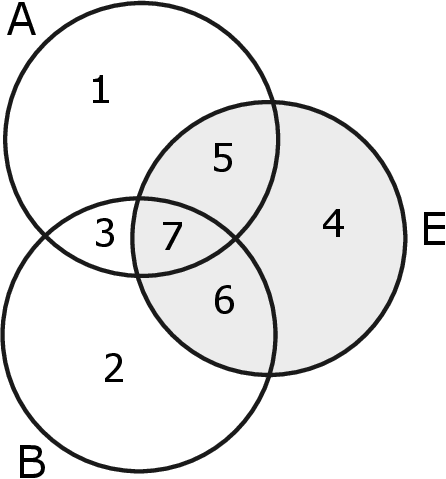}}
\end{center}

Since $B$ is measurable, $[4567]=[45]+[67]$.
Since $A$ is measurable, we have $[67]=[6]+[7]$.
Using these equations and sub-additivity, we have
$$[4567] \leq [456]+[7] \leq [45]+[6]+[7]=[45]+[67]=[4567].$$
So, we have equalities everywhere and we have proved what we want.
\endproof

We note in particular that when $A$ and $B$ are measurable
and $A \subset B$ we have $\mu(B \cap A) + \mu(B -A)=\mu(B)$, or
$\mu(B-A)=\mu(B)-\mu(A)$.   Also, by induction, finite unions and intersections of
measurable sets are measurable.

\begin{lemma}
  \label{countable}
  The countable union of measurable sets is measurable.
\end{lemma}

\startproof
Suppose $S=\bigcup S_n$ is a countable union of
measurable sets.  Since finite unions of measurable sets are
measurable, we can suppose that we have an
increasing union.  Now we define
$T_n=S_n-S_{n-1}$ and $T_1=S_1$.
Since these sets are all measurable we have
$\mu(T_n)=\mu(S_n)-\mu(S_{n-1})$.  Hence
\begin{equation}
  \label{finitesum}
\sum_{k=1}^n \mu(T_k)=\mu(S_n) \leq \mu(S) \leq \sum_{k=1}^{\infty}
\mu(T_k).
\end{equation}
The first inequality is monotonicity and the last is subadditivity.
Letting $n \to \infty$ and using the fact that a bounded monotone
sequence converges,
\begin{equation}
  \label{finitesum2}
  \sum_{k=1}^{\infty} \mu(T_k) =\lim \mu(S_n) \leq \mu(S) \leq
  \sum_{k=1}^{\infty} \mu(T_k).
\end{equation}
But then, in fact, we have equality everywhere in Equation \ref{finitesum2}.
In particular the series converges.  Hence we can choose $n$ so large that
\begin{equation}
  \label{small}
  \mu(S_n^c-S^c)=\mu(S-S_n) =\mu\bigg(\bigcup_{k=n}^{\infty} T_k\bigg) \leq \sum_{k=n}^{\infty}
  \mu(T_k)<\epsilon,
\end{equation}
Combining Equation \ref{small} with the measurability of $S_n$ and
sub-additivity,
$$\mu(E) \leq \mu(E \cap S) + \mu(E \cap S^c) \leq
\mu(E \cap S_n) + \mu(E \cap S_n^c) + $$ $$\mu(S-S_n) + \mu(S_n^c-S^c) \leq
\mu(E)+2\epsilon.$$
Since $\epsilon$ is arbitrary, we have
$\mu(E)=\mu(E \cap S) + \mu(E \cap S^c)$.
\endproof

\noindent
{\bf Proof of Statement 1:\/}
Any open set is a countable union of dyadic cubes.
So, Lemma \ref{cube} and Lemma \ref{countable} say that open 
sets are measurable.  Hence closed sets are measurable.
This fact combines with Lemma \ref{countable}
to show that all Borel sets are measurable.
\endproof

\noindent
{\bf Proof of Statement 2:\/}
This follows from Equations \ref{finitesum} and \ref{finitesum2}.
\endproof

\noindent
{\bf Proof of Statement 3:\/}
We get open $U \supset S$ with $\mu(U)-\mu(S)<\epsilon/2$ by
taking a suitable cover of $S$ by open cubes and then taking their
union.
We get a compact $K \subset S$ by applying the same
construction to $S^c$ and
then taking complements.  The measurability of all sets involved
guarantees that $\mu(S)-\mu(K)<\epsilon/2$.
Combining these, and again using measurability, we have
$\mu(U-K)=\mu(U)-\mu(K)<\epsilon$.
\endproof

\subsection{Baby Fubini}

Here we prove a weak special case of the Fubini-Tonelli Theorem
[{\bf F\/}, p. 65] which is
sufficient for Theorem \ref{null2}.  We identify the set of horizontal
lines of $[0,1]^2$ with $[0,1]$.

\begin{lemma}
  \label{FUBINI}
  Suppose $S$ is measurable.
  Let $F_{t,S}$ denote the set of
  horizontal lines $L$ such that $\mu(L \cap S)>t$.
  If $\mu(S)<t^2$ then $\mu(F_{t,S}) \leq t$.
\end{lemma}

\startproof
In view of Statement 3 of Theorem \ref{borel}, and the
corresponding result for null sets, it suffices to prove
this result when $S$ is open.
If $S$ is a finite union of dyadic squares, we subdivide
and reduce to the case when they all have the same size and
belong to a common grid of $[0,1]^2$.   In this case, we get the
result by counting.
In general $S=\bigcup Q_n$ is a carpet.
Also, the intersection
$L \cap S$ is a carpet
for each horizontal line $L$.
Let $Q^n$ be the union
of the first $n$ squares.
By the finite case,
$\mu(F_{t,Q^n})  \leq t$.
By Statement 2 of Theorem \ref{borel},
$\mu(\bigcup F_{t,Q^n})\leq t$ and
$\mu(L \cap S)=\lim \mu(L \cap Q^n)$ for
each horizontal line $L$.
Hence $F_{t,S} \subset \bigcup F_{t,Q^n}$.
Hence $\mu(F_{t,S}) \leq t$.
\endproof

\subsection{Covering Results}

Our proofs here take their ideas from  [{\bf F\/}, \S 3.4].
We work with closed intervals.
A {\it Besicovich cover\/} $\cal B$ of $S \subset [0,1]$ is a union of
intervals, such that each point of $S$ is centered on some
interval of $\cal B$ and each interval of $\cal B$ is centered
on some point of $S$.  Let $|I|$ denote the length
of an interval $I$.

\begin{lemma}
  \label{cover1}
  Suppose $\cal B$ is a Besicovich cover of $S$.
  Then there exists a subset of $\cal B$ consisting
  of disjoint intervals whose total length-sum is
    at least $\mu(S)/3$.
\end{lemma}

\startproof
Greedily choose intervals,
always picking a largest one that is disjoint
from the previous ones picked. Let
$\{I_j\}$ be this collection.  Let $S' = \bigcup I_j$.
Let $3I_j$ be the interval obtained by dilating
$I_j$ by a factor of $3$ about its midpoint.
Any $x \in S-S'$ is the center of some $J$ of $\cal B$ not
picked by the algorithm.   But then $J$ intersects
some $I_j$ with $|I_j| \geq |J|$.
This forces $x \in 3I_j$.  Hence $\{3I_j\}$
covers $S$.  Hence $\sum |3I_j| \geq \mu(S)$.
Hence $\sum |I_j| \geq \mu(S)/3$.
\endproof

An interval $J$ is $\delta$-{\it porous\/} if
$\mu(J \cap K) <(1-\delta)|J|$.  A point
$p \in [0,1]$ is $\delta$-{\it porous\/} if $p$ is the center of
arbitrarily small $\delta$-porous intervals.
Finally, $K$ is {\it porous\/} if for some $\delta>0$
every point of $K$ is $\delta$-porous.

\begin{lemma}[Porous]
  If $S$ is measurable and porous then $\mu(S)=0$.
\end{lemma}

\startproof
Suppose $\mu(S)>0$. By Statement 3 of Theorem \ref{sandwich}
we can find $K \subset S \subset U$
with $K$ compact, $U$ open, and
$\mu(U-K)<\epsilon$.  If we take $\epsilon$ small,
we have $\mu(K)>0$.
The set $K$ remains porous.

Since $K$ is compact, there is some $\lambda>0$ such that every interval
of length $\lambda$ centered at a point of
$K$ lies in $U$.
We take a Besicovich
covering $\cal B$ of $K$ by $\delta$-porous
intervals all of length less than $\lambda$.
Letting $\{I_j\}$ be as in
Lemma \ref{cover1}, we have
$\delta_j:=\mu(I_j \cap U)- \mu(I_j \cap K) \geq |I_j| -
(1-\delta)|I_j|  =\delta |I_j|$.
From this calculation, we see that the set $U-K$ intersects
$\bigcup I_j$ in a set of size at least $\sum \delta_j \geq
\delta \mu(K)/3$.  This is a contradiction for small enough
$\epsilon$.
\endproof

We say that a covering
$\cal B$ by intervals is {\it renewable\/} if,
for every $\epsilon>0$, every
point of $S$ is the (left or right)
endpoint of an interval in $\cal B$
having length less than $\epsilon$. 
Let $S \Delta T=(S-T) \cup (T-S)$.

\begin{theorem}[Vitali]
  \label{symm}
  Suppose  $\cal B$ is a renewable cover for a 
 Borel set $S$.   For any $\epsilon>0$
  there is a disjoint collection of intervals
  $\{I_j\}$ of $\cal B$ such that
  $\mu(S \Delta T)<\epsilon$ where
  $T=\bigcup I_j$.
\end{theorem}

\startproof
If $S$ is null, the result is trivially true.
So, assume $\mu(S)>0$.
By Statement 3 of Theorem \ref{sandwich} we have
$K \subset S \subset U$ where $K$ is compact and $U$ is open
and $\mu(U-K)<\epsilon/2$.
There is some $\lambda>0$ so that every interval of
length $\lambda$ having an endpoint in $K$ belongs to $U$.
We make $\{I_j\}$ by length-greedily picking disjoint
intervals of $\cal B$, having endpoints in $K$ and length
less than $\lambda$.
Let $T=\bigcup I_j$.   By monotonicity,
$\mu(T-S)<\epsilon/2$.

Let $p \in K'=K-T$.
Let $I$ be an interval in
$\cal B$ having length less than $\lambda$ and
$p$ as an endpoint.
Let $J$ be the interval
centered at $p$ having length $8|I|$.
The interval $I$ was not picked in the
algorithm, so there is some interval
$I_j$ with $|I_j| \geq |I|$, such that
$p$ is at most $|I_j|$ away from an endpoint
  of $I_j$.
But then $\mu(I_j \cap J) \geq |I| = |J|/8$.
Hence $\mu(K' \cap J) \leq (7/8) |J|$.
Hence $p$ is $(1/8)$-porous with respect to $K'$.
Hence $K'$ is $(1/8)$-porous.
Also, $K'$ is a Borel set and hence measurable.
We conclude from the Porous Lemma that $\mu(K')=0$.
Hence
$  \mu(S-T) \leq \mu(S-K)+ \mu(K')<\epsilon/2$.
This combines with $\mu(T-S)<\epsilon/2$ to give
$\mu(S \Delta T)<\epsilon$.
\endproof

\section{Proofs of the Analytic Results}

\subsection{Proof of Lemma \ref{nice}}
\label{niceproof}

Let $\phi: [0,1]^2 \to \R$ be a continuous map.  We extend
$\phi$ to be  continuous on $\R^2$. 
  Let $D \subset (0,1)^2$ be the set where $\partial
  \phi/\partial x$ exists.  Let
  \begin{equation}
  \label{diffquot}
  \phi_{\epsilon}(x,y)=\frac{\phi(x+\epsilon,y)-\phi(x,y)}{\epsilon}
\end{equation}
Let's call $\phi$ {\it rationally differentiable\/} at $(x,y)$
if $\{\phi_{\epsilon}(x,y)\}$ is a Cauchy sequence relative to any
sequence of rational $\epsilon$ with   $\epsilon \to 0$.  
Because $\phi$ is continuous, $\phi$ is differentiable at
$(x,y)$  if and only if $\phi$ is rationally differentiable
at $(x,y)$.

In other terms, this is true
iff for each $n \geq 1$ there is an
integer $m \geq 1$ such that if $0<|a|,|b|<1/m$ are rational
then $(x,y) \in U_{m,a,b,n}$, the set such 
that $|\phi_a(x,y)-\phi_b(x,y)|<1/n$.
Since $\phi_a$ and $\phi_b$ are
continuous, $U_{m,a,b,n}$ is an open set, and hence a Borel set.
So, $D$ is a Borel set because
$$D=\bigcap_{n=1}^{\infty} \bigcup_{m=1}^{\infty} U_{m,n}, \hskip 30 pt
U_{m,n}=\bigcap_{0<|a|,|b|<1/m} U_{m,a,b,n}.$$

\subsection{Proof of Theorem \ref{null2}}
\label{fullproof}

For any set $Y$ we let
$Y^c$ denote its complement in the relevant domain.
Given $S \subset [0,1]^2$ we
let $F_S$ and $F_{t,S}$ respectively denote the set of
horizontal lines $L$ such that
$\mu(L \cap S)>0$ and $\mu(L \cap S)>t$.

Suppose $\mu(F_S)>0$.
Then there is some $t>0$ such that $\mu(F_{t,S})>t$.
If $\mu(S)=0$ then also $\mu(S)<t^2$ and we
contradict Lemma \ref{FUBINI}.

Suppose $S$ is a Borel set and $\mu(S)>0$.
Then $\mu(S^c)<t^2$ for some $t<1$.
By Lemma \ref{FUBINI}, we have
$\mu(F_{t,S^c}) \leq t$.
But if $L \notin F_{t,S^c}$ then $\mu(L \cap S^c)\le t$, and hence
$\mu(L \cap S)\ge 1-t$.  So
$(F_{t,S^c})^c \subset F_{1-t,S}$.
Therefore
$\mu(F_{1-t,S}) \geq 1-t$.  Since
$F_{1-t,S} \subset F_S$, we have $\mu(F_S)>0$.
Hence, if $\mu(F_S)=0$ then $\mu(S)=0$.
      
\subsection{Proof of Theorem \ref{stiff}}
\label{stiffproof}

Let $S$ be the stretchy set for
$A: {\cal J\/} \to (0,\infty)$.
If $S$ is fat then
$\mu(S)>1/N$ for some $N$.  Let
$\cal B$ be a Besicovich cover of $S$
by $N^2$-stretched intervals.  Let
$\{J_j\}$ be the set of disjoint intervals produced by
Lemma \ref{cover1}.   For large $N$ we have a
contradiction:
$A([0,1]) \geq \sum A(J_j) \geq N^2 \mu(S)/3 \geq N/3.$
So, $S$ is null.

\subsection{Proof of Theorem \ref{AC}}
\label{ACproof}

\noindent
{\bf Proof Outline:\/}
First we give the proof modulo four
auxiliary lemmas.
Lemma \ref{leb2} below shows that
an AC function is the difference
between two monotone AC functions.
Lemma \ref{leb} below shows that
a monotone AC function is differentiable
except on a null set.
These two results combine to show that
an arbitrary AC function is differentiable
away from a null set.

Let $f$ be an AC function and
let $A$ be the set where $f'=0$.   
Let $B=[0,1]-A$. Suppose $B$ is null.
Lemma \ref{CO3} below says in particular that
the image of a null set under an
AC map is null.  Hence
$f(B)$ is null.
By Lemma \ref{CO4} below, the set
$f(A)$ is null. But then
$f([0,1])=f(A) \cup f(B)$ is both null and connected.
Hence $f([0,1])$ is a single point. In particular, $f(0)=f(1)$.
Hence, if $f(0) \not = f(1)$ then $B$ is fat.
\newline

Now we will take care of the details of these lemmas. 
For convenience we repeat the definition of an
AC function.
  Suppose that $f: [0,1] \to \R$ is a continuous function.
Let $I=\{I_1,...,I_n\}$ denote a finite list of
intervals of $[0,1]$ having pairwise disjoint
interiors.  We call $I$ a
{\it partial partition\/}.  Let $|I|=\sum |I_k|$.
We define $I_k'$ to be the interval bounded
by the two points of $f(\partial I_k)$.
We define $I'=\{I_1',...,I_n'\}$ and
$|I'|=\sum |I_k'|$. 
The function $f$ is {\it AC\/} (absolutely continuous)
if, for
each $\epsilon>0$, there is some
$\delta>0$ such that
$|I|<\delta$ implies that $|I'|<\epsilon$.

\begin{lemma} 
  \label{CO3}
  Let $f$ be an AC function.
  Let $S \subset [0,1]$.
For each $\epsilon>0$ there is some
$\delta>0$ such that $\mu(S)<\delta$ implies that
$\mu(f(S))<\epsilon$. Hence, if $S$ is null then $f(S)$ is null.
\end{lemma}

\startproof 
Suppose $\mu(S)<\delta$.  We can find an open set $U$
such that $S \subset U$ and $\mu(U)<\delta$.
Note that $U$ is a countable union of intervals with disjoint interiors.
Let $U^n$ denote the union
of the first $n$ intervals of $U$.
We construct a partial partition $I=I^n$ as follows.
  For each compact connected component $C$ of ${\rm closure\/}(U^n)$,
we include in $I^n$ an interval connecting
points of $C$ where $f$ respectively achieves its min and max.
By construction, $|I|<\delta$ and
$f(U^n) \subset V^n:=\bigcup I'_k$.
Choosing $\delta$ small enough, we have
$|I'|<\epsilon$.  Hence $\mu(V^n)<\epsilon$.
Also by construction $V^n \subset V^{n+1}$ for all $n$.
Since $V = \bigcup V^n$ is an increasing union of
sets, each a finite union of intervals,
$\mu(V) \leq \epsilon$.  
But $f(S) \subset V$.  Hence $\mu(f(S))<\epsilon$.
\endproof

Our next lemma is usually stated for functions of
{\it bounded variation\/}.  In the notation of the
lemma below, this means that $v(f,[0,1])<\infty$.
Here we have
stronger hypotheses and a stronger conclusion.

\begin{lemma}
  \label{leb2}
  If $f$ is AC then $f=f_+-f_-$ where $f_{\pm}$ is monotone and AC.
  \end{lemma}

\startproof
Given an interval $Y \subset [0,1]$ we define
the {\it variation\/} $v(f,Y)$ to be the
supremum of $|I'|$ taken over all partial partitions $I$ of $Y$.

We claim first that $v(f,[0,1])<\infty$.
If $v(f,[0,1])=\infty$ then for any $\delta>0$ we can find
an interval $Y$ such that $|Y|<\delta$ and $v(f,Y) \geq 1$.
This contradicts the fact that $f$ is AC. Hence
$v(f,[0,1])<\infty$.

Given the monotonicity properties of the variation, we
conclude that $v(f,Y)<\infty$ for every interval $Y
\subset [0,1]$. We also note an additivity property:
For $0 \leq a<b$ we have
\begin{equation}
  \label{add}
  v(f,[0,b])=v(f,[0,a])+v(f,[a,b]).
  \end{equation}
Let $f_+(x)=v(f,[0,x])$.  By construction,
$f_+$ is increasing.  By Equation \ref{add},
$$f_+(b)-f_+(a)=v(f,[a,b]) \geq f(b)-f(a).$$
Hence $f_-=f_+-f$ is also increasing.

It remains to show that $f_+$ and $f_-$ are
AC. Since $f$ is AC, and sums and differences
of AC functions are AC, it suffices to prove
that $f_+$ is AC.
Let us reformulate the
definition of AC first.  Instead of defining
$|I'_k|$ as the distance between the two points of
$f(\partial I_k)$ we might
have defined $|I'_k|=v(f,I_k)$.   Since the AC criterion
is already defined in terms of taking a supremum
over all partial partitions, the inclusion of the
further supremum in the definition of $|I'_k|$
changes nothing.  This alternate
definition picks out the same class of functions
as AC.  But this definition is formulated entirely
in terms of the variation function.
Now, by Equation \ref{add}, 
$$v(f_+,[a,b])=f_+(b)-f_+(a)=v(f,[a,b]).$$
Since $f$ and $f_+$ have the same variation
function and $f$ is AC, so is $f_+$.
\endproof

The next result is the workhorse in our proof
of Theorem \ref{AC}.  It
actually holds for all monotone functions.
See [{\bf Fa\/}].  The case of monotone AC functions
is much easier, though the core idea is the same as
in the general case.

\begin{lemma}
  \label{leb}
  If $f$ is  monotone and AC then $f'$ exists outside a null set.
    \end{lemma}

    \startproof
    Adding a linear function to $f$, we reduce to the case
        when $f$ is strictly monotone.
Let $E_{a,b}$ denote the set of points $p \in [0,1]$ such
that
\begin{enumerate}
\item For all $n \in \N$ the point $p$ is the endpoint of
  an interval $I$ such that $\mu(I)<1/n$ and
  $\mu(f(I))<a \mu(I)$.
\item For all $n \in \N$ the point $p$ is the endpoint of
  an interval $I$ such that $\mu(I)<1/n$ and
  $\mu(f(I))>b \mu(I)$.
\end{enumerate}
To prove this lemma, it suffices to prove
that $\mu(E_{a,b})=0$ for all $0 \leq a<b$.
Suppose some $S=E_{a,b}$ has $\mu(S)>0$.
Since $f$ is continuous, $S=E_{a,b}$ is defined by
a countable collection of open conditions. Hence $S$ is a Borel set.

Let $\cal B$ be the renewable cover of $S$ made from
the intervals in Item 1.
For any $\delta>0$ let
$\{I_j\}$ and $T=\bigcup I_j$
be as in Theorem \ref{symm} so that
$\mu(S \Delta T)<\delta$.
By construction $\mu(f(T)) \leq a \mu(T)$.
By Lemma \ref{CO3}, we can make
$\mu(S \Delta T)$ and
$\mu(f(S \Delta T))$ as small
as we like by shrinking $\delta$.  Hence
$\mu(S) \leq a \mu(S)$.

Running the same argument with Item 2 in place of Item 1, we get
$\mu(S) \geq b \mu(S)$.  Since
$b \mu(S) \leq a \mu(S)$ and
$a<b$ we have $\mu(S)=0$.
\endproof

Our final result does not need the AC hypothesis.  

\begin{lemma}
  \label{CO4}
  Suppose $f: [0,1] \to \R$ is continuous.  Suppose also that
  $f'$ exists and equals $0$ on a set $A \subset [0,1]$. Then
  $f(A)$ is null.
\end{lemma}

\startproof
Each $x \in A$ is the midpoint of an
interval $I_x$ such that
$|f(I_x)|<\epsilon |I_x|$.
Let $\cal B$ be the Besicovich cover made from these
intervals.  We introduce an auxiliary
Besicovich cover
${\cal B\/}^*$ of $A$, obtained by
shrinking all the intervals in ${\cal B\/}$ by a
factor of $3$ about their midpoints.  Let
$\{I_j^*\}$ be the disjoint collection
of intervals of ${\cal B\/}^*$ produced by
the greedy algorithm.   Our proof of Lemma
\ref{cover1} shows that the corresponding union of
dilated intervals $\{I_j\}$ is a cover of
$A$.  By construction, $\sum |I_j| \leq 3$.
But now we can say $\mu(f(A))<3\epsilon$. Since $\epsilon$
is arbitrary, $f(A)$ is null.  
\endproof

\section{Discussion}
\label{scope}

Usually one defines a compact hyperbolic manifold
to be a compact Riemannian manifold that is locally
isometric to $\H^3$.  By the Cartan-Hadamard Theorem,
this more typical definition agrees with the definition
in \S \ref{backg}.

Our proof showed that the extension map
$h$ lies in $\I$.
Hence the isometry from $M_1$ to $M_2$ we
get induces the same isomorphism from
$\pi_1(M_1)$ to $\pi_1(M_2)$ that the BL map
$f$ does.

If we just assume that $f: M_1 \to M_2$ is a
homotopy equivalence, then the lifted map
$H$ is still a quasi-isometry. Theorem \ref{extend}
still works in this case. So two homotopy
equivalent compact hyperbolic $3$-manifolds
are isometric.

The same proof works for compact
hyperbolic $n$-manifolds when $n \geq 4$.
The key change is in
Theorem \ref{heart}.   When $n=4$ we would
let $\cal J$ denote the set of axis-aligned squares of $[0,1]^2$
instead of the set of intervals of $[0,1]$
and then Equation \ref{lv} would be in terms of
volume rather than area.

The proof also works for
finite volume hyperbolic $n$-manifolds, when we have
$n \geq 3$.
In the finite volume case, you have
to avoid zooming into 
cusps.

\newpage

\section{References}

\noindent
[{\bf F\/}] G. B. Folland,
{\it Real Analysis: Modern Techniques and Their Applications\/},
John Wiley \& Sons, New York (1984).
\vspace{8pt}
\newline
\noindent
[{\bf Fa\/}] C. Faure,
{\it An elementary proof of the differentiability almost everywhere
of monotone functions\/},
Real Analysis Exchange {\bf 15} (1989–90)
\vspace{7 pt}
\newline
\noindent
[{\bf GS\/}] S. Gouëzel and V. Shchur,
{\it A corrected quantitative version of the Morse lemma\/},
Journal of Functional Analysis {\bf 277}, pp.\ 1258–1268 (2019).
\vspace{8pt}
\newline 
\noindent
[{\bf LV\/}] O. Lehto and K. I. Virtanen,
{\it Quasiconformal Mappings in the Plane\/},
Second Edition, Springer--Verlag, New York (1973).
\vspace{8pt}
\newline 
\noindent
[{\bf M\/}] G. D. Mostow,
{\it Quasi-conformal mappings in $n$-space and the rigidity of hyperbolic space forms\/},
Inst. Hautes \'{E}tudes Sci. Publ. Math. No. 34, pp. 53--104 (1968).
\vspace{8pt}
\newline 
\noindent
[{\bf S\/}] R. E. Schwartz,
{\it The quasi-isometry classification of rank one lattices\/},
Inst. Hautes \'{E}tudes Sci. Publ. Math. No. 82, pp. 133--168 (1995).
\vspace{8pt}
\newline 
\noindent
[{\bf T\/}] W. P. Thurston,
{\it The Geometry and Topology of Three-Manifolds\/},
Princeton University Lecture Notes (1978).
\vspace{8pt}
\newline
[{\bf Z\/}] L. Zaj\'\i\v{c}ek,
{\it Fr\'echet differentiability via partial Fr\'echet differentiability\/},
Comment. Math. Univ. Carolin. {\bf 64}, No. 2, pp. 185--207 (2023).

\end{document}